\documentclass[a4paper,10pt,reqno]{amsart}
\usepackage[a4paper,tmargin=40mm,bmargin=35mm,hmargin=30mm]{geometry}
\usepackage[T1]{fontenc}
\usepackage{lmodern}
\usepackage{amsmath}
\usepackage{amssymb}
\usepackage{amsthm}
\usepackage{stmaryrd}
\usepackage{tikz}
\usepackage{booktabs}
\usepackage{multirow}
\usepackage{multicol}
\usepackage{here}
\usepackage{aliascnt}
\usepackage{hyperref}
\usepackage[noabbrev]{cleveref}

\newcommand{\NewTheorem}[2]{
	\newaliascnt{#1}{TheoremEnvironment}
	\newtheorem{#1}[#1]{#1}
	\aliascntresetthe{#1}
	\crefname{#1}{#1}{#2}
	\Crefname{#1}{#1}{#2}
}

\theoremstyle{definition}

\NewTheorem{Definition}{Definitions}
\NewTheorem{Remark}{Remarks}
\NewTheorem{Axiom}{Axioms}
\NewTheorem{Example}{Examples}
\NewTheorem{Observation}{Observations}
\NewTheorem{Convention}{Conventions}
\NewTheorem{Notation}{Notations}
\NewTheorem{Setting}{Settings}
\NewTheorem{Question}{Questions}
\NewTheorem{Answer}{Answers}
\NewTheorem{Conjecture}{Conjectures}
\NewTheorem{Problem}{Problems}
\NewTheorem{Solution}{Solutions}
\NewTheorem{Goal}{Goals}
\NewTheorem{Comment}{Comments}
\NewTheorem{Aim}{Aims}
\NewTheorem{Caution}{Cautions}
\NewTheorem{Exercise}{Exercises}

\theoremstyle{plain}
\NewTheorem{Proposition}{Propositions}
\NewTheorem{Lemma}{Lemmas}
\NewTheorem{Theorem}{Theorems}
\NewTheorem{Corollary}{Corollaries}

\crefname{enumi}{}{}
\Crefname{enumi}{}{}
\creflabelformat{enumi}{(#2#1#3)}
\crefname{enumii}{}{}
\Crefname{enumii}{}{}
\creflabelformat{enumii}{(#2#1#3)}
\crefname{enumiii}{}{}
\Crefname{enumiii}{}{}
\creflabelformat{enumiii}{(#2#1#3)}

\makeatletter
\renewcommand{\p@enumii}{}
\renewcommand{\p@enumiii}{}
\makeatother

\numberwithin{equation}{section}
\crefname{equation}{}{}
\Crefname{equation}{}{}
\creflabelformat{equation}{(#2#1#3)}

\newcommand{\SwapSymbols}[1]{
	\expandafter\let\expandafter\temporarysymbol\csname #1\endcsname
	\expandafter\let\csname #1\expandafter\endcsname\csname var#1\endcsname
	\expandafter\let\csname var#1\endcsname\temporarysymbol
}

\SwapSymbols{epsilon}
\SwapSymbols{phi}
\SwapSymbols{Gamma}
\SwapSymbols{Delta}
\SwapSymbols{Theta}
\SwapSymbols{Lambda}
\SwapSymbols{Xi}
\SwapSymbols{Pi}
\SwapSymbols{Sigma}
\SwapSymbols{Upsilon}
\SwapSymbols{Phi}
\SwapSymbols{Psi}
\SwapSymbols{Omega}

\newcommand{\bF}{\mathbf{F}}

\newcommand{\CA}{\mathcal{A}}

\newcommand{\C}{\mathcal{C}}
\newcommand{\DD}{\mathcal{D}}

\newcommand{\cF}{\mathcal{F}}

\newcommand{\CI}{\mathcal{I}}

\newcommand{\cL}{\mathcal{L}}
\newcommand{\cM}{\mathcal{M}}

\newcommand{\cS}{\mathcal{S}}
\newcommand{\cT}{\mathcal{T}}

\newcommand{\CX}{\mathcal{X}}

\newcommand{\To}{\rightarrow}

\DeclareMathOperator{\Hom}{Hom}

\DeclareMathOperator{\Ext}{Ext}

\DeclareMathOperator{\Modd}{Mod}

\DeclareMathOperator{\modd}{mod}
\DeclareMathOperator{\omodd}{\underline{mod}}

\DeclareMathOperator{\pd}{pd}

\DeclareMathOperator{\Ker}{Ker}
\DeclareMathOperator{\Coker}{Coker}

\DeclareMathOperator{\im}{Im}

\title{Torsion theory of coherent functors}
\subjclass[2010]{18E15 (Primary), 16P60, 16D90, (Secondary)}
\keywords{functors category, coherent functor, torsion theory, radical functor}

\author{Mohammad Khazaei and Reza Sazeedeh}
\address{Department of Mathematics, Urmia University, P.O.Box: 165, Urmia, Iran}
\email{rsazeedeh@ipm.ir}
\address{Department of Mathematics, Urmia University, P.O.Box: 165, Urmia, Iran}
\email{mohammad.khazaei@gmail.com}

\begin{document}


\begin{abstract}
Let $\C$ be an additive category with cokernels and let $\Modd(\C)$ be the category of additive functors from $\C^{op}$ to the category Ab of abelian groups. Let $\modd(\C)$ be the full subcategory of $\Modd(\C)$ consisting of coherent functors. In this paper, we first study some basic properties of pseudo-kernel of morphisms in $\C$. When $\C$ has pseudo-kernels, $\modd(\C)$ is abelian and then, in this case, we study radical functors, half exact functors, left exact functors and injective objects in $\modd(\C)$. At last, we extend the results for $\Modd(\C)$.   
\end{abstract}
\maketitle

\section{introduction}
The notion of coherent first time in 1940 was introduced by Cartan [C] and Oka [O] on several complex variables. Serre in 1955 in his famous paper [S], showed that the same notions could be carried over to algebraic geometry and since then coherent sheaves and their cohomology have been ubiquitous in algebraic geometry. In 1962, Auslander [A], developed the notion coherent in functors category. Hartshorne [H], explained the general theory of coherent functors on the category of finitely generated modules over a nothereian ring, in order to study coherent sheaves on projective space. The various characterization of functors category was given by Krause [K1]. He investigated certain subcategories which are defined in terms of coherent functors. 

Throughout this paper we assume that $\C$ is an additive category with cokernels. We denote by $\Modd(\C)$ the category of additive functors $F:\C^{{\rm op}}\To{\rm Ab}$ which is an abelian category. A functor $F:\C^{{\rm op}}\To{\rm Ab}$ in $\Modd(C)$ is called {\it coherent} if there exists an exact sequence $\Hom_{\C}(-,A)\To \Hom_{\C}(-,B)\To F\To 0$ of functors such that $A$ and $B$ are objects in $\C$. We denote by $\modd(\C)$, the full category of $\Modd(\C)$ which consists coherent functors. By Yoneda lemma for a functor $F$ with the presentation $\Hom_{\C}(-,A)\To \Hom_{\C}(-,B)\To F\To 0$, there is a correspondence morphism $A\To B$ and if we set $v(F)=\Coker (A\To B)$, then $v:\modd(\C)\To \C$ is an additive functor and according to [K2, Lemma 2.8], if $\C$ has cokernels, then $v$ has a right adjoint functor $h_{\C}$ by the assignment $h_{\C}(C)=\Hom_{\C}(-,C)$. We denote by $\omodd(\C)$ a full subcategory of $\modd(\C)$ consisting of all functors $F$ such that $v(F)=0$.

 In section 2, we find a characterization of functors in $\omodd(\C)$ and we show that $(\omodd(\C),\DD)$ is a torsion theory in which $\DD$ is a full subcategory of $\modd(C)$ consisting of all functors which send epimorphisms to monomorphisms (cf. \cref{epimono} and \cref{tor}). We recall that a morphism $B\To A$ is a pseudo-kernel of a morphism $A\To C$ in $\C$ provided that the induced sequence of functors $\Hom_{\C}(-,B)\To \Hom_{\C}(-,A)\To \Hom_{\C}(-,C)$ is exact. We study some basic properties
of pseudo-kernels which have a key role in the other sections. We show that $v$ preserves pseudo-kernels if and only if $\C$ is abelian (cf. \cref{ab}).

In section 3, we assume that $\C$ has pseudo-kernels and so by [F, Theorem 1.4], $\modd(\C)$ is an abelian category. A full subcategory $\CX$ of $\modd(\C)$ is called Serre if it is closed under subobject, quotients and extension. In this section we assume that $\omodd(\C)$ is a Serre subcategory of $\modd(\C)$. We show that a functor $G\in\modd(\C)$ is left exact if and only if $\Ext_{\modd(\C)}^i(F,G)=0$ for $i=0,1$ and all $F\in\omodd(\C)$ (cf. \cref{leftex}). We prove that for any functor $F\in\modd(\C)$, there exist a unique left exact functor $L(F)$ and a morphism $\eta_F:F\To L(F)$ such that $\Ker\eta_F$ and $\Coker\eta_F$ are in $\omodd(\C)$ and so we define a radical functor on $\modd(\C)$ with $r(F)=\Ker\eta_F$ for any $F\in\modd(\C)$ (cf. \cref{local}). We show that $r:\modd(\C)\To\modd(\C)$ is left exact and $r:\modd(\C)\to\omodd(C)$ is a right adjoint functor of the inclusion functor. We denote by $\cF_r$ a subclass of $\modd(\C)$ consisting of all functors $F$ such that $r(F)=0$ and we show that $\DD=\cF_r$ and if $\C$ has kernels, every $F\in\cF_r$ has projective dimension $\leq 1$ (cf. \cref{ttt}). Finally, in \cref{p}, we prove that $r(F)=0$ if and only if $F$ has a presentation $\Hom_{\C}(-,A)\stackrel{\theta}\To \Hom_{\C}(-,B)\stackrel{\pi}\To F\To 0$ such that $v(\theta)$ is pseudo-kernel of $v(\pi)$ for any functor $F\in\modd(C)$. 

In section 4, we study half exact functors in $\modd(\C)$. A functor $F$ in $\modd(\C)$ is said to be half exact if for any exact sequence $ A\To B\To C\To 0$ of objects of $\C$ with $A\To B$ a pseudo-kernel of $B\To C$, the seuence $F(C)\To F(B)\To F(A)$ is exact. We show that a functor $G$ is half exact if and only if $\Ext^1_{\modd(\C)}(F,G)=1$ for all $F\in\omodd(\C)$. In particular if $\omodd(\C))$ is Serre and $G$ is half exact, then so is $r(G)$. For any functors $F$ and $G$ in $\omodd(\C)$, there is a homomorphism $\theta_{F,G}:F(v(G))\To\Hom_{\modd(\C)}(F,G)$ which is natural in $F$ and $G$ and we show that $F\in\DD$ if and only if $\theta_{F,G}$ is monic for all $G\in\omodd(C)$ (cf. \cref{gy} and \cref{c1}).  We also show that $F$ is half exact if and only if $\theta_{F,G}$ is epic for all $G\in\modd(\C)$; in particular, $\theta_{F,G}$ is isomorphism for all $G\in\modd(\C)$ if and only if $F$ is left exact (cf. \cref{c2}). In the rest of this section we study injective objects in $\modd(\C)$. In \cref{rinj}, we show that if $G\in\modd(\C)$ is injective, then so is $r(G)$. We also show that if $\C$ is abelian and $A$ is an object of $\C$ such that ${\rm id}A=n$, then $\Ext^n_{\C}(-,A)$ is injective in $\modd(\C)$. Moreover, if $n\geq 1$, then $\Ext^n_{\C}(-,A)$ is in $\omodd(\C)$ (cf. \cref{inje}).

For any subcategory $\cS$ of $\modd(\C)$, we denote by $\overrightarrow{\cS}$ the full subcategory of $\Modd(\C)$ which consists of direct limit $\underrightarrow{\rm lim} F_i$ with each $F_i\in\cS$. In section 5, we show that $(\overrightarrow{\omodd(\C)},\overrightarrow{\DD})$ is a hereditary torsion theory of $\Modd(\C)$ (cf. \cref{torm}). In \cref{st}, we prove that the subcategory of $\overrightarrow{\omodd(\C)}$-closed objects in $\Modd(\C)$ is precisely $\overrightarrow{\cL}$, where $\cL$ is a subcategory of $\modd(\C)$ consisting of left exact functors. If $\C$ has pseudo-kernels, in \cref{lefex}, we show that $\overrightarrow{\cL}$ is the class of left exact functors in $\Modd(\C)$. For any functor $F\in\Modd(\C)$, we prove that there exists a morphism $\eta_F:F\To L(F)$ such that $L(F)\in\overrightarrow{\cL}$ and $\Ker\eta_F$ and $\Coker\eta_F$ are in $\overrightarrow{\omodd(\C)}$ (cf. \cref{ex}). This leads us to define a radical functor $r:\Modd(\C)\To \Modd(\C)$ which coincides with $r:\modd(\C)\to\modd(\C)$ on $\modd(C)$. At the end of this paper, we show that if $\C$ has direct limits, then every $F\in\Modd(\C)$ has $\cL$-preenvelopes.
	
\section{categories with cokernels}
Let $\mathcal{U}$ be a universe in the sense of Grothendieck [G]. We say that a category $\C$ is a $\mathcal{U}$-category if the objects and the morphisms of $\C$ are both sets which are elements of $\mathcal{U}$ (for details see [A] and [G]). Throughout this paper $\C$ is an additive $\mathcal{U}$-category and we denote by $\Modd(\C)$ the category of additive functors from $\C^{\rm op}$ to the category Ab of abelian groups. For any objects $A, B$ in $\C$, the set of all morphisms from $A$ to $B$ is denoted $\Hom_{\C}(A,B)$ which is an abelian group. It follows from Yoneda lemma that $\Hom_{\C}(-,B)$ is a projective object of $\Modd(\C)$ for any object $B$ in $\C$.

Let $\C$ be an additive category. A functor $F$ in $\Modd(\C)$ is called {\it coherent} if there exist objects $A$ and $B$ in $\C$ such that there is an exact sequence of functors $\Hom_{\C}(-,A)\To \Hom_{\C}(-,B)\To F\To 0$. We denote by $\modd(\C)$ a full subcategory of $\Modd(\C)$ consisting of coherent functors.

Using Yoneda's lemma, there exists a full and faithful faunctor $h_{\C}:\modd(\C)\To\C, A\mapsto \Hom_{\C}(-,A)$ and according to [K2, Lemma 2.8], if $\C$ has cokernels, then $h_{\C}$ has a left adjoint functor $v:\modd(\C)\To \C$. It follows from [K2, Universal Property 2.1] that $v=1^{\star}_{\C}$ is right exact and $v h_{\C}=1_{\C}$.  To be more precise, if $F\in\modd(\C)$, then there is a presentation $\Hom(-,X)\To\Hom(-,Y)\To F\To 0$ and we define $v(F)=\Coker(X\To Y)$. We also define a full subcategory $\omodd(\C)$ of $\modd(\C)$ which is $\omodd(\C)=\{F\in\modd(\C)| \hspace{0.1cm} v(F)=0\}$. It is straightforward to show that $\omodd(\C)$ is closed under quotients and extensions. Throughout this section $C$ has cokernels. We start this section by a results about  $\omodd(\C)$.

\medskip
	\begin{Proposition}\label{epimono}
	Let $F\in\modd(\C)$. Then the following are equivalent.
	
	{\rm(i)} $F\in\omodd(\C)$.
	
	{\rm(ii)} there exists an epimorphism $A_2\To A_1$ such that $(-,A_2)\To (-,A_1)\To F\To 0$ is exact in  $\modd(\C)$.
	
	{\rm(iii)} $\Hom_{\modd(\C)}(F,G)=0$ for those $G\in$ $\modd(\C)$ which send epimorphisms to monomorphisms.
		\end{Proposition}
	\begin{proof}
	$\rm{(i)}\Rightarrow\rm{(ii)}$ and $\rm{(ii)}\Rightarrow\rm{(iii)}$ are clear. $\rm{(iii)}\Rightarrow\rm{(i)}$ Since $v$ is a left adjoint of $h$ and by the assumption we have $\Hom_{\C}(v(F),v(F))\cong \Hom_{\modd(\C)}(F, h_{\C}v(F))=0$ which implies that $v(F)=0$.
				\end{proof}

	\medskip

		Let $\DD$ be a subcategory of $\modd(\C)$ consisting of all functors which send epimorphisms to monomorphisms. Then we have the following proposition.
	
	\begin{Proposition}\label{tor}
	 $(\omodd(\C),\DD)$ is a torsion theory.
	\end{Proposition}
	\begin{proof}
	Assume that $G\in\modd(\C)$	such that $\Hom_{\modd(C)}(F,G)=0 $ for all functors $F\in\omodd(\C)$. For any epimorphism $A\To B$ in $\C$, there exists $F\in\omodd(\C)$ such that the following sequence is exact $$\Hom_{\C}(-,A)\to \Hom_{\C}(-,B)\To F\To 0.$$ According to Yoneda lemma, $G(B)\to G(A)$ is monomorphism and so the proof completes by using \cref{epimono}.
	\end{proof}

	\medskip

	We recall  that a morphism $B\To A$ is a {\it pseudo-kernel} of a morphism $A\To C$ in $\C$ provided that the induced sequence  of functors $\Hom_{\C}(-,B)\To \Hom_{\C}(-,A)\To \Hom_{\C}(-,C)$ is exact. This concept introduced by Freyd and he proved that $\C$ has pseudo-kernels iff $\modd(\C)$ is abelian (cf. [F, Theorem 1.4]). The {\it pseudo-cokernel} is defined dually. In the rest of this section we assume that $\C$ has pseudo-kernels.  We study some basic properties for this category.

\medskip

\begin{Proposition}\label{uniq}
Cokernel of pseudo-kernels of any morphism of $\C$ is unique up to isomorphisms. 
\end{Proposition}
\begin{proof}
 Assume that $\omega _1:\Omega_1\To A$ and $\omega_2:\Omega_2\To A$ are pseudo-kernels of $\gamma:A\To B$. By the definition there exist morphisms $\theta_1:\Omega_1\To\Omega_2$ and $\theta_2:\Omega_2\to\Omega_1$ such that $\omega_1\theta_2=\omega_2$ and $\omega_2\theta_1=\omega_1$. Assume that $\beta_i:A\To C_i$ are cokernels of $\omega_i$ for $i=1,2$. Since $\beta_1\omega_2=\beta_1\omega_1\theta_2=0$ and $\beta_2\omega_1=\beta_2\omega_2\theta_1=0$, there exist the morphisms $\alpha_1:C_1\To C_2$ and $\alpha_2:C_2\To C_1$ such that $\alpha_1\alpha_2=1_{C_2}$ and $\alpha_2\alpha_1=1_{C_1}$.
\end{proof}

\medskip

\begin{Lemma}\label{psco}
If $\omega:\Omega\To A$ is a pseudo-kernel of a morphism, then $\omega$ is a pseudo-kernel of $\Coker\omega$.
\end{Lemma}
\begin{proof}
Assume that $\omega$ is pseudo-kernel of $\gamma:A\To B$ and $\alpha:A\To C$ is the cokernel of $\omega$. Then there exists a morphism $\theta:C\To B$ such that $\theta\alpha=\gamma$.
Assume that $\delta:X\To A$ is any morphism such that $\alpha\delta=0$. Then $\gamma\delta=0$; and hence there is a morphism $\eta:W\To X$ such that $\delta\eta=\omega$. 
\end{proof}

\medskip

\begin{Lemma}
		If any morphism of $\C$ factors through an epimorphism and a monomorphism and $\C$ has pseudo-kernels, then it has kernels.
		\end{Lemma}
		\begin{proof}
		Assume that $\phi:Y\To Z$ is a morphism and $\psi:X\To Y$ is a pseudo-kernel of $\phi$. Then there is a factorization of $\psi$ as
		$X\stackrel{\psi'} \To I\stackrel{i}\To Y$ where $i$ is monomorphism. It is straightforward to show that $i$ is a pseudo-kernel of $\phi$ and since $i$ is monomorphism, it is kernel of $\phi$.
	\end{proof}

\medskip
The following theorem shows that when a category having cokernels and pseudo-kernels is abelian. 

\begin{Theorem}\label{ab} 
Let $\C$ have pseudo-kernels. Then the following conditions are equivalent:

${\rm (i)}$  The functor $v$ preserves pseudo-kernels.

${\rm (ii)}$ The functor $v:\modd(\C)\To \C$ is left exact.

${\rm (iii)}$ $\C$ is abelian 
\end{Theorem}
	\begin{proof}
		(i)$\Leftrightarrow (ii)$ follows from [K2, Lemma 2.6].
	(ii)$\Rightarrow$ (iii). Assume that $A\To B$ is a morphism in $\C$ and assume that $F$ is the kernel of $\theta:\Hom_{\C}(-,A)\To\Hom_{\C}(-,B)$. Since $v$ is left exact, it preserves monomorphisms and pseudo-kernels; and hence $v(F)$ is the kernel of $v(\theta)=A\To B$. On the other hand, since $\modd(\C)$ is abelian, we have $\Coker\Ker\theta\cong\Ker\Coker\theta$ and since $v$ is exact, we have $\Coker\Ker v(\theta)\cong\Ker\Coker v(\theta)$. Further, since $\modd(\C)$ is abelian, it has finite products, and since $v$ is exact and an additive functors, it is easy to prove that $\C$ has finite product. (iii)$\Leftrightarrow$(ii)	follows from [A].
			\end{proof}

	\medskip

	Let $f:\C\To \DD$ be a functor between additive categories. We denote by $f^{\star}:\modd(\C)\to\modd(\DD)$ the unique right exact functor extending $f$. By [K2, Universal Property 2.1] such a functor always exists. 

\medskip

\begin{Proposition}
Let $\C$ and $\DD$ be categories with cokernels and $f:\C\To\DD$ be an additive functor.  Then $f$ preserves epimorphisms iff $f^{\star}(\Ker(v_{\C}))\subseteq\Ker(v_{\DD})$. Moreover, if $\C$ has kernels, $\DD$ has pseudo-kernels and $f$ preserves pseudo-kernels, then $f$ preserves epimorphisms iff $f$ is exact.
\end{Proposition}
\begin{proof}
Assume that $f$ preserves epimorphisms and $F\in\Ker(v_{\C})$. Then there exists an exact sequence 
$\Hom_{\C}(-,B)\To \Hom_{\C}(-,A)\To F\To 0$ such that $B\To A$ is an epimorphism in $\C$. In view of the assumption 
$f(B)\To f(A)$ is epimorphism in $\DD$ so that $f^{\star}(F)=\Coker(\Hom_{\DD}(-,f(B))\To \Hom_{\DD}(-,f(A))\in\Ker(v_{\DD})$. The proof of the converse is similar. To prove the second claim, using [K2, Lemma 2.6], $f$ is left exact and so it suffices to show that it is right exact. Given an exact sequence $0\To B\To A\To C\To 0$, by the assumption $f(B)\To f(A)$ is pseudo-kernel of $f(A)\To f(C)$ and hence $\Hom_{\DD}(-, f(A))\To \Hom_{\DD}(-,f(B))\To \Hom_{\DD}(-,f(C))$ is exact. On the other hand, $F=\Coker(\Hom_{\C}(-,A)\To \Hom_{\C}(-,C))\in\Ker(v_{\C})$ and so using the first assertion, $f(B)\To f(A)\To f(C)\To 0$ is exact. The converse is clear. 
\end{proof}


		\section{radical functors}
		
		 A full subcategory $\CX$ of $\modd(\C)$ is said to be {\it Serre} if it is closed under quotients, subobjects and extensions. Throughout this section we assume that $\C$ has pseudo-kernels and cokernels and $\omodd(\C)$ is a Serre subcategory of $\modd(\C)$

			\begin{Lemma}\label{copse}
	Let $B\To A$ be an epimorphim of $\C$. If $C\To B$ is a pseudo-kernel of $B\To A$, then $C\To B\To A\To 0$ is exact (i.e $B\To A=\Coker(C\To B)$). 
	\end{Lemma}
	\begin{proof}
	According to \cref{epimono}, there exists an exact sequence $\Hom_{\C}(-,C)\To \Hom_{\C}(-,B)\To \Hom_{\C}(-,A)\To F\To 0$ of functors such that $F\in\omodd(\C)$. If $B\to D$ is the cokernel of $C\To B$, then $B\To A$ factors through an epimorphism $D\To A$.  It follows from \cref{psco} that $C\To B$ pseudo-kernel of $B\To D$ and so there exists $F_0\in\omodd(\C)$ such that the following diagram with exact rows commuts
	
	\begin{center}\setlength{\unitlength}{3cm}
	\begin{picture}(5.3,.5)
\put(.5,0.45){\makebox(0,0){$\Hom_{\C}(-,C)$}}
\put(1.6,0.45){\makebox(0,0){$\Hom_{\C}(-,B)$}}
 \put(2.7,0.45){\makebox(0,0){$\Hom_{\C}(-,D)$}}
\put(3.7,0.45){\makebox(0,0){$F_0$}}
\put(4.3,0.45){\makebox(0,0){$0$}}

\put(.35,.38){\line(0,-1){.135}}
\put(.32,.38){\line(0,-1){.135}}

\put(1.45,.38){\line(0,-1){.135}}
\put(1.42,.38){\line(0,-1){.135}}

\put(2.55,.38){\vector(0,-1){.135}}
\put(2.58,.3){$\gamma$}

\put(3.7,.39){\vector(0,-1){.135}}
\put(3.73,.3){$\theta$}

\put(.85,.45){\vector(1,0){0.37}}
\put(.85,.2){\vector(1,0){0.37}}
\put(1.95,.45){\vector(1,0){0.37}}
\put(1.95,.2){\vector(1,0){0.37}}
\put(3.05,.45){\vector(1,0){0.5}}
\put(3.05,.2){\vector(1,0){0.5}}
\put(3.785,.45){\vector(1,0){0.45}}
\put(3.785,.2){\vector(1,0){0.45}}

\put(.5,0.2){\makebox(0,0){$\Hom_{\C}(-,C)$}}
\put(1.6,0.2){\makebox(0,0){$\Hom_{\C}(-,B)$}}
 \put(2.7,0.2){\makebox(0,0){$\Hom_{\C}(-,A)$}}
\put(3.7,0.2){\makebox(0,0){$F$}}
\put(4.3,0.2){\makebox(0,0){$0.$}}
\end{picture}
\end{center}
		Hence there is an isomorphism $\Ker\gamma\cong\Ker\theta$ and since $\omodd(\C)$ is Serre, $\Ker\theta\in\omodd(\C)$. Hence \cref{tor} implies that $\Ker\gamma=0$ so that the epimorphism $D\to A$ is isomorphism.  
		\end{proof}
	
	\medskip
	
	We now find a characterization of left exact functors.
	\begin{Proposition}\label{leftex}	
	Let $G\in\modd(\C)$. Then $\Ext^i_{\modd(\C)}(F,G)=0$ for $i=0,1$ and all $F\in\omodd(\C)$ if and only if $G$ is left exact.  
	\end{Proposition}
	\begin{proof}
		Assume that $\Ext^i_{\modd(\C)}(F,G)=0$ for $i=0,1$ and all $F\in\omodd(\C)$. There exists a morphism $A_1\To A_0$ and an exact sequence of functors $\Hom_{\C}(-,A_1)\To\Hom_{\C}(-,A_0)\To G\To 0$. Assume that $A_2\To A_1$ is a pseudo-kernel of $A_1\to A_0$ and also assume that $B=\Coker(A_2\To A_1)$. Then $v(G)=\Coker(B\To A_0)$ and by virtue of \cref{psco}, the morphism $A_2\To A_1$ is a pseudo-kernel of $A_1\to B$. Therefore, we have the following commutative diagram with the exact rows 

\begin{center}\setlength{\unitlength}{3cm}
\begin{picture}(5.3,.5)

\put(.5,0.45){\makebox(0,0){$\Hom_{\C}(-,A_2)$}}
\put(1.6,0.45){\makebox(0,0){$\Hom_{\C}(-,A_1)$}}
 \put(2.7,0.45){\makebox(0,0){$\Hom_{\C}(-,B)$}}
\put(3.7,0.45){\makebox(0,0){$G_0$}}
\put(4.3,0.45){\makebox(0,0){$0$}}

\put(.35,.38){\line(0,-1){.135}}
\put(.32,.38){\line(0,-1){.135}}

\put(1.45,.38){\line(0,-1){.135}}
\put(1.42,.38){\line(0,-1){.135}}

\put(2.55,.38){\vector(0,-1){.135}}
\put(2.58,.3){$\gamma$}
\put(2.55,.13){\vector(0,-1){.145}}
\put(3.7,.39){\vector(0,-1){.135}}
\put(3.73,.3){$\theta$}

\put(.85,.45){\vector(1,0){0.37}}
\put(.85,.2){\vector(1,0){0.37}}
\put(1.95,.45){\vector(1,0){0.37}}
\put(1.95,.2){\vector(1,0){0.37}}
\put(3.05,.45){\vector(1,0){0.5}}
\put(3.05,.2){\vector(1,0){0.5}}
\put(3.785,.45){\vector(1,0){0.45}}
\put(3.785,.2){\vector(1,0){0.45}}

\put(.5,0.2){\makebox(0,0){$\Hom_{\C}(-,A_2)$}}
\put(1.6,0.2){\makebox(0,0){$\Hom_{\C}(-,A_1)$}}
 \put(2.7,0.2){\makebox(0,0){$\Hom_{\C}(-,A_0)$}}
\put(3.7,0.2){\makebox(0,0){$G$}}
\put(4.3,0.2){\makebox(0,0){$0$}}
\put(2.7,-.05){\makebox(0,0){$\Hom_{\C}(-,v(G))$}}
\end{picture}
\end{center}

 where by using \cref{epimono} we have $G_0\in\omodd(\C)$. Since $\omodd(\C)$ is Serre, we have $\Ker\theta\in\omodd(\C)$ and an isomorphism $\Ker\gamma\cong\Ker \theta$. Then it follows from \cref{epimono} that $\Ker\theta=\Ker\gamma=0$. We observe that $\Hom_{\C}(-, B)\To\Hom_{\C}(-,A_0)$ is monic so that the induced morphism $B\To A_0$ is monic. If we put $F=\Coker(\Hom_{\C}(-,B)\To \Hom_{\C}(-,A_0)$, then $v(G)=v(F)=\Coker(B\To A_0)$ so that $0\To \Hom_{\C}(-,B)\To \Hom_{\C}(-,A_2)\To \Hom_{\C}(-,v(F))$ is exact. Then using the above diagram, we have an exact sequence $0\To G_0\To G\To \Hom_{\C}(-,v(G))\To G_1\To 0$ where $G_1=\Coker(\Hom_{\C}(-,A_0)\To\Hom_{\C}(-,v(G)))$ and so using \cref{epimono}, it belongs to $\omodd(\C)$. Since $G_0\in\omodd(\C)$, the assumption implies that $G_0=0$ and $0\To G\To \Hom_{\C}(-,v(G))\To G_1\To 0$ splits so that $G_1=0$. Conversely, assume that $G$ is left exact and $F\in\omodd(\C)$. Then there exist objects $B$ and $C$ and an epimorphism $\alpha:B\To C$ such that
	$\Hom_{\C}(-,B)\To \Hom_{\C}(-,C)\To F\To 0$ is exact. If $K$ is the pseudo-kernel of $\alpha$, it follows from \cref{copse} that $K\To B\To C\To 0$ is exact and so there is an exact sequence of functors $\Hom_{\C}(-,K)\To \Hom_{\C}(-,B)\To \Hom_{\C}(-,C)\To F\To 0$. Now, application of the functor $\Hom_{\C}(-,G)$ and using the assumption and Yoneda lemma, the result follows. 
	 \end{proof}
	
	\medskip
	\begin{Proposition}\label{local}
	For any $F\in\modd(\C)$, there exists a left exact functor $L(F)$ and a morphism $\eta_F:F\To L(F)$ such that $\Ker\eta_F$ and $\Coker\eta_F$ are in $\omodd(\C)$. Furthermore, $L(F)$ and $\eta_F$ with the mentioned property is unique up to isomorphism.  
	\end{Proposition}
	\begin{proof}
	By a similar proof of \cref{leftex}, there exists an exact sequence of functors $E(F):0\To F_0\To F\stackrel{\eta_F}\To \Hom_{\C}(-,v(F))\To F_1\To 0$ such that $F_0$ and $F_1$ are in $\omodd(\C)$. Then we consider $L(F)=h_{\C}v(F)=\Hom_{\C}(-,v(F))$. A similar proof of [A, Proposition 3.4], shows that the above exact sequence is unique up to isomorphism.  
	\end{proof}
	 
	\medskip

A preradical functor $r$ of $\modd(\C)$ is a subfunctor of $1_{\modd(\C)}$. The preradical $r$ is called idempotent if  $r^2=r$ and it is called radical if for any functor $F$ in $\modd(\C)$ we have $r(F/r(F))=0$. If $\omodd(\C)$ is Serre, then we define a preradical of $\modd(\C)$ as $r(F):=\Ker\eta_F$ for any $F$ in $\modd(\C)$. We observe that $\eta_{(-)}:(-)\To L(-)$ is a natural transform. To be more precise, if $f:F\To G$ is a morphism of functors in $\modd(\C)$, then applying the functor $\Hom_{\C}(-,L(G))$ to $E(F)$, there exist unique morphisms $L(f):L(F)\To L(G)$ and $r(F)\To r(G)$ such that the following diagram is commutative

\begin{center}\setlength{\unitlength}{3cm}
\begin{picture}(5.3,.5)

\put(.5,0.4){\makebox(0,0){$0$}}
\put(1.2,0.4){\makebox(0,0){$r(F)$}}
 \put(2,0.4){\makebox(0,0){$F$}}
\put(3,0.4){\makebox(0,0){$L(F)$}}
\put(4,0.4){\makebox(0,0){$\Coker\eta_F$}}
\put(4.6,0.4){\makebox(0,0){$0$}}

\put(1.05,0.18){\makebox(0,0){$r(f)$}}
\put(1.2,.35){\vector(0,-1){.28}}
\put(1.9,0.18){\makebox(0,0){$f$}}
\put(2,.35){\vector(0,-1){.28}}
\put(2.85,0.18){\makebox(0,0){$L(f)$}}
\put(3,.35){\vector(0,-1){.28}}
\put(4,.35){\vector(0,-1){.28}}

\put(.6,.4){\vector(1,0){0.45}}
\put(1.37,.4){\vector(1,0){0.57}}
\put(2.4,0.45){\makebox(0,0){$\eta_F$}}
\put(2.4,0.05){\makebox(0,0){$\eta_G$}}
\put(2.1,.4){\vector(1,0){0.757}}
\put(3.17,.4){\vector(1,0){0.57}}
\put(4.25,.4){\vector(1,0){0.3}}

\put(.6,.4){\vector(1,0){0.45}}
\put(1.37,.4){\vector(1,0){0.57}}
\put(2.1,.4){\vector(1,0){0.757}}
\put(3.17,.4){\vector(1,0){0.57}}
\put(4.25,.4){\vector(1,0){0.3}}

\put(.6,.0){\vector(1,0){0.45}}
\put(1.37,.0){\vector(1,0){0.57}}
\put(2.1,.0){\vector(1,0){0.757}}
\put(3.17,.0){\vector(1,0){0.57}}
\put(4.25,.0){\vector(1,0){0.3}}

\put(.5,0){\makebox(0,0){$0$}}
\put(1.2,0){\makebox(0,0){$r(G)$}}
 \put(2,0){\makebox(0,0){$G$}}
\put(3,0){\makebox(0,0){$L(G)$}}
\put(4,0){\makebox(0,0){$\Coker\eta_G$}}
\put(4.6,0){\makebox(0,0){$0$}}
\end{picture}
\end{center}

Clearly, $r$ is an idempotent preradical functor. We also define the {\it torsion-free class} $\mathcal{F}_r$ a subclass of $\modd(\C)$ consisting of those $F\in\modd(\C)$ which $\eta_F$ is monic; or equivalently $r(F)=0$.

	\medskip

	\begin{Lemma}\label{adje}
	The functor $r:\modd(\C)\To\omodd(\C)$ is a right adjoint of the inclusion functor $i:\omodd(\C)\To \modd(\C)$
	\end{Lemma}
	\begin{proof}
	For any $G\in\omodd(\C)$ and $F\in\modd(\C)$, applying $\Hom_{\modd(\C)}(G,-)$ to the  exact sequence
	$$0\to r(F)\To F\stackrel{\eta_F}\To h_{\C}v(F)\To\Coker{\eta_F}\To 0$$
	we deduce that $\Hom_{\omodd(\C)}(G,r(F))=\Hom_{\modd(\C)}(G,r(F))\cong\Hom_{\modd(\C)}(G,F)$ so that $r:\modd(\C)\To\omodd(\C)$ is the right adjoin functor of the inclusion functor $i:\omodd(\C)\To\modd(\C)$. 
	\end{proof}
	
	Since $r$ is a right adjoint of $i$, the functor $r:\modd(\C)\To \omodd(\C)$ is left exact and since $h_{\C}$ is a right adjoint of $v$,  $r(F)$ is the largest subfunctor of $F$ contained in $\omodd(\C)$.

	\medskip
	\begin{Lemma}\label{left}
The functor $r:\modd(\C)\To\modd(\C)$ is left exact. In particular, $r$ is radical.
\end{Lemma}
\begin{proof}
Let $0\To F_1\stackrel{\iota}\To F_2\stackrel{\theta}\To F_3\To 0$ is an exact sequence in $\modd(\C)$ and so there is an sequence of functors $0\To r(F_1)\stackrel{r(\iota)}\To r(F_2)\stackrel{r(\theta)}\To r(F_3)$. We notice that $r(\iota)$ is monic and if $K$ is the kernel of $r(\theta)$, since $\omodd(\C)$ is Serre, we have $K\in\omodd(C)$ and there is a monomorphism $\alpha:r(F_1)\To K$; and also there is a monomorphism $h:K\To F_1$ such that $h\circ\alpha=i_{F_1}:r(F_1)\To F_1$ is the inclusion functor. Since $r(F_1)$ is the largest subfunctor of $F_1$ contained in $\omodd(\C)$, $\alpha$ is an isomorphism. In order to prove the second assertion, applying the functor $r$ to the exact sequence $0\To F/r(F)\To h_{\C}v(F)$ yields that $r(F/r(F))=0$.    
\end{proof}

\medskip

The following theorem shows that the subcalss $\cF_r$ coincides with the class $\DD$ mentioned in Section 2.

\begin{Theorem}\label{ttt}
Let $F$ be a functor in $\C$ and consider the following conditions.\\
$\rm (i)$ $F\in\mathcal{F}_r$.\\
$\rm (ii)$ $F\in\mathcal{D}$.\\
$\rm (iii)$ $\pd F\leq 1$.
Then $\rm(i)$ and $\rm(ii)$ are equivalent. The implication $\rm(iii)\Rightarrow\rm(i)$ holds and if $\C$ has kernels, then the implication $\rm(iii)\Rightarrow\rm(i)$ holds too.
\end{Theorem}
\begin{proof}
(i)$\Rightarrow$ (ii). There exists an exact sequence of functors 
$0\To F\stackrel{\eta_F}\To (-,v(F))$ in $\modd(\C)$. For any exact sequence of objects $A\stackrel{f}\To B\To 0$ in $\C$, we have 
$\Hom(f,v(F)) \eta_F(B)=\eta_F(A) F(f)$ which implies that $F(f)$ is monic. (ii)$\Rightarrow$ (i). According to \cref{tor}, we have $\Hom_{\C}(r(F), F)=0$ and so the exact sequence $0\To r(F)\To F\To \Hom_{\C}(-,v(F))$ implies that $r(F)=0$. (iii)$\Rightarrow$ (i). If pd$F=0$, then $F$ is a direct summand of $\Hom_{\C}(-,A)$ where $A$ is an object of $\C$. It follows from \cref{epimono} that $(r(\Hom_{\C}(-,A)),\Hom_{\C}(-,A))=0$ which implies that $r(\Hom_{\C}(-,A))=0$; and hence $r(F)=0$. If pd$F= 1$, then there exists an exact sequence of functors $0\To \Hom_{\C}(-,A)\To \Hom_{\C}(-,B)\To F\To 0$ which yields an exact sequence $0\To A\To B\To v(F)\To 0$ of object of $\C$. Application of $\Hom_{\C}(-,)$ induces the following exact sequence of functors 
$$0\To \Hom_{\C}(-,A)\To \Hom_{\C}(-,B)\To \Hom_{\C}(-,v(F))$$ which induces an exact sequence $0\To F\To \Hom_{\C}(-,v(F))$ of functors in $\C$. Application of the functor $\Hom_{\C}(r(F),-)$ and using \cref{epimono} yield $r(F)=0$. (i)$\Rightarrow$ (iii). Since $r(F)=0$, there is the exact sequence $0\To F\stackrel{\eta_F}\To\Hom_{\C}(-,v(F)$ of functors and since $F\in\modd(\C)$, there is an exact sequence of functors
$\Hom_{\C}(-,A)\stackrel{\alpha}\To \Hom_{\C}(-,B)\stackrel{\pi}\To F\To 0$. Assume that $f:C\To B$ is the kernel the morphism $B\To v(F)$. Then there exists a unique morphism $g:A\to C$ such that $f\circ g=v(\alpha)$. Further, there is an exact sequence $0\To\Hom_{\C}(-,C)\stackrel{\gamma}\To \Hom_{\C}(-,B)\To\Hom_{\C}(-,v(F))$ in $\modd(\C)$ where $v(\gamma)=f$. We have the following commutative diagram with exact rows where $v(\beta)=g$ \\\\ 

\begin{center}\setlength{\unitlength}{3cm}
\begin{picture}(4.5,.7)
\put(0,0){\makebox(0,0){$0$}}
 \put(0.85,0){\makebox(0,0){$\Hom_{\C}(-,C)$}}
\put(2.1,0){\makebox(0,0){$\Hom_{\C}(-,B))$}}
 \put(3.5,0){\makebox(0,0){$\Hom_{\C}(-,v(F))$}}

\put(0.1,0){\vector(1,0){0.38}}
\put(1.25,0){\vector(1,0){0.42}}
 \put(2.52,0){\vector(1,0){0.514}}

 \put(0.9,0.7){\makebox(0,0){$\Hom_{\C}(-,A)$}}
\put(2.1,0.7){\makebox(0,0){$\Hom_{\C}(-,B)$}}
 \put(3.55,0.7){\makebox(0,0){$F$}}
\put(4.3,0.7){\makebox(0,0){$0$}}

\put(1.3,.75){$\alpha$}
\put(1.3,.05){$\gamma$}
\put(1.31,.7){\vector(1,0){0.4}}
 \put(2.52,.7){\vector(1,0){0.95}}
\put(3.75,.7){\vector(1,0){.4}}

 \put(0.85,.6){\vector(0,-1){0.452}}
 \put(2.1,.6){\vector(0,-1){.452}}
\put(3.52,.6){\vector(0,-1){0.452}}
\put(.87,.315){$\beta$}
 \put(2.17,0.315){$1_{\Hom_{\C}(-,B)}$}
\put(3.57,0.315){$\eta_F$}
\end{picture}
\end{center}
Since, $\eta_F$ is monic, $\beta$ is epic; and hence; we deduce that $F=\Coker\alpha=\Coker\gamma$ which forces $\pd F\leq 1$. 
\end{proof}

\medskip
We show that the class $\cF_r$ can be specified in terms of pseudo-kernels. We find a characterization of functors in $\cF_{r}$. 
	\begin{Theorem}\label{p}
		Let $F\in\modd(\C)$. Then $r(F)=0$ if and only if there exists a presentation $\Hom_{\C}(-,A)\stackrel{\theta}\To\Hom_{\C}(-,B)\To F\To 0$ of $F$ such that $v(\theta)$ is the pseudo-kernel of $\Coker v(\theta)$. 
		\end{Theorem}
		\begin{proof}
		Assume that $r(F)=0$. Then there exists a presentation $\Hom_{\C}(-,A)\stackrel{\theta}\To\Hom_{\C}(-,B)\stackrel{\pi}\To F\To 0$ of $F$ and an exact sequence of functors $\Hom_{\C}(-,A)\stackrel{\theta}\To \Hom_{\C}(-,B)\stackrel{f}\To \Hom_{\C}(-,v(F))$ such that $f=\eta_{F}\circ\pi$ and $\eta_F$ is monic. Therefore $v(\theta):A\To B$ is the pseudo-kernel of $v(f):B\To v(F)$. It now follows from \cref{psco} that $v(\theta)$ is pseudo-kernel of $v(\pi)$. Conversely assume that $\Hom_{\C}(-,A)\stackrel{\theta}\To\Hom_{\C}(-,B)\stackrel{\pi}\To F\To 0$ is a presentation of $F$ such that $v(\theta)$ is a pseudo-kernel of $v(\pi)$. Then we have an exact sequence $\Hom_{\C}(-,A)\stackrel{\theta}\To \Hom_{\C}(-,B)\stackrel{f}\To \Hom_{\C}(-,v(F))$ of functors and we set $I=\im \theta$. Now, there exists the following commutative diagram with the exact rows  \\\\
		
		\begin{center}\setlength{\unitlength}{3cm}
\begin{picture}(4.5,.7)
\put(0,0){\makebox(0,0){$0$}}
 \put(0.85,0){\makebox(0,0){$I$}}
\put(2.1,0){\makebox(0,0){$\Hom_{\C}(-,B)$}}
 \put(3.5,0){\makebox(0,0){$\Hom_{\C}(-,v(F))$}}

\put(0.1,0){\vector(1,0){0.6}}
\put(0.95,0){\vector(1,0){0.75}}
 \put(2.475,0){\vector(1,0){0.55}}

 \put(0.9,0.7){\makebox(0,0){$\Hom_{\C}(-,A)$}}
\put(2.1,0.7){\makebox(0,0){$\Hom_{\C}(-,B)$}}
 \put(3.55,0.7){\makebox(0,0){$F$}}
\put(4.1,0.7){\makebox(0,0){$0$}}

\put(1.3,.75){$\theta$}

\put(1.28,.7){\vector(1,0){0.4}}
 \put(2.48,.7){\vector(1,0){0.95}}
\put(3.65,.7){\vector(1,0){.4}}

 \put(0.85,.6){\vector(0,-1){0.452}}
 \put(2.1,.6){\vector(0,-1){.452}}
\put(3.52,.6){\vector(0,-1){0.452}}

\put(.92,0.315){$p$}
 \put(2.17,0.315){$1_{\Hom_{\C}(-,B)}$}
\put(3.57,0.315){$\alpha$}
\end{picture}
\end{center}
				
	which implies that $\alpha$ is monic. Since $r$ is left exact, we have $r(F)=0$.
	\end{proof}

	\medskip

	\section{half exact functors}
	Throughout this section, we assume that $\C$ has pseudo-kernels and cokernels. 
			
 A functor $F$ in $\modd(\C)$ is called {\it half exact} if for any exact sequence $ A\To B\To C\To 0$ of objects of $\C$ with $A\To B$ a pseudo-kernel of $B\To C$, the sequence $F(C)\To F(B)\To F(A)$ is exact. 

In the following proposition, we show that half exact functors coincides with those which have been mentioned in the abelian categories.

		\begin{Proposition}\label{quasihalf}
Assume that $\C$ has kernels and $F$ is a functor in $\modd(\C)$. Then $F$ is half exact if for any exact sequence $0\To A\To B\To C\To 0$, $F(C)\To F(B)\To F(A)$ is exact. 
\end{Proposition}
\begin{proof}
Assume that $A\To B\To C\To 0$ is an exact sequence such that $A\To B$ is a pseudo-kernel of $B\To C$. If $K\To B$ is the kernel of $B\To C$, it follows from \cref{uniq} that $0\To K\To B\To C\To 0$ is an exact sequence and $K$ is a direct summand of $A$. Using the assumption, one can easily show that $F(C)\To F(B)\To F(A)$ is exact.  
\end{proof}

	\medskip
	We now give a characterization of half exact functors.
	\begin{Lemma}\label{1}
	Let $G\in\modd(\C)$. Then $G$ is half exact if and only if $\Ext^1_{\modd(\C)}(F,G)=0$ for all $F\in\omodd(\C)$. In particular, if $\omodd(\C)$ is Serre and $G$ is half exact, then so is $r(G)$.
	\end{Lemma}
	\begin{proof}
	For any $F\in\omodd(\C)$ there exists an epimorphism $B\To C$ such that the sequence $\Hom_{\C}(-,B)\To\Hom_{\C}(-,C)\To F\To 0$ is exact. Then there exists an object $A$ such that $A\To B$ is a pseudo-kernel of $B\To C$. It follows from \cref{copse} that $A\To B\To C\To 0$ is exact so that $\Hom_{\C}(-,A)\To\Hom_{\C}(-,B)\To\Hom_{\C}(-,C)\To F\To 0$. Now, the result follows by applying $\Hom_{\modd(\C)}(-,G)$ and using the Yoneda lemma. To prove the second assertion, for any $F\in\omodd(\C)$, applying the functor $\Hom_{\modd(\C)}(F,-)$ to the exact sequence $0\To r(G)\To G\To G/r(G)\To 0$ and using \cref{ttt}, \cref{tor}, the assumption and the first assertion, the result follows. 
	\end{proof}

	 \medskip
	\begin{Proposition}\label{2}
	Let $\omodd(\C)$ be a Serre subcategory of $\modd(\C)$ and $G\in\modd(\C)$. Then the following conditions hold.
	
	${\rm(i)}$ $G$ preserves epimorphisms to monomorphisms if and only if $\eta_G$ is monic.
	
	${\rm (ii)}$  $G$ is left exact if and only if $\eta_G$ is isomorphism.
	\end{Proposition}
	\begin{proof}
	(i) is clear by \cref{tor} and \cref{ttt}.(ii) If $G$ is left exact, then according to \cref{epimono}, we have $\Hom_{\modd(\C)}(r(G),G)=0$ and so $r(G)=0$. Then using \cref{leftex}, the exact sequence $0\To G\stackrel{\eta_G}\To h_{\C}v(G)\To F\To 0$ splits so that $F=0$. The converse of the implication is clear.
		\end{proof}

In the following theorem we give a generalization of the Yoneda lemma.
 \medskip
\begin{Theorem}\label{gy}
Let $F$ and $G$ be functors in $\modd(\C)$. Then there is a homomorphism $\theta_{F,G}:F(v(G))\To \Hom
_{\C}(G,F)$ of abelian groups which is natural in $F$ and $G$. Moreover, if $F$ is in $\DD$, then $\theta_{F,G}$ is monic.  
\end{Theorem}
\begin{proof}
Since $G\in\modd(\C)$, there exists an exact sequence of functors $\Hom_{\C}(-,B)\To \Hom_{\C}(-,A)\To G\To 0$. Application of the functor $v(-)$ yields the exact sequence $B\stackrel{\alpha}\To A\stackrel{\pi}\To v(G)\To 0$ of objects of $\C$. If we put $H=\im(\Hom_{\C}(-,B)\To \Hom_{\C}(-,A))$ and $C=v(H)$, then we have a natural epimorphism 
$\beta:B\To C$ and a morphism $\gamma:C\To A$ such that $\alpha=\gamma\beta$ and $ C\stackrel{\gamma}\To A\stackrel{\pi}\To v(G)\To 0$ is an exact sequence of objects in $\C$. Then, there is a sequence $F(v(G))\stackrel{F(\pi)}\To F(A)\stackrel{F(\gamma)}\To F(C)$ of abelian groups and further, by Yoneda lemma, there exist an exact sequence $0\To \Hom_{\C}(G,F)\To F(A)\stackrel{F(\alpha)}\To F(B)$ of abelian groups and a unique homomorphism $\theta_{F,G}: F(v(G))\To (G,F)$ such that the following diagram is commutative\\
 
\begin{center}\setlength{\unitlength}{3cm}
\begin{picture}(4.5,.7)
\put(0,0){\makebox(0,0){$0$}}
 \put(0.85,0){\makebox(0,0){$\Hom_{\modd(\C)}(G,F))$}}
\put(2.1,0){\makebox(0,0){$F(A)$}}
 \put(3.5,0){\makebox(0,0){$F(B)$}}

\put(0.08,0){\vector(1,0){0.293}}
\put(1.35,0){\vector(1,0){0.582}}
 \put(2.34,0){\vector(1,0){0.98}}

 \put(0.9,0.7){\makebox(0,0){$F(v(G))$}}
\put(2.1,0.7){\makebox(0,0){$F(A)$}}
 \put(3.55,0.7){\makebox(0,0){$F(C)$}}

\put(1.2,.7){\vector(1,0){0.7}}
 \put(2.4,.7){\vector(1,0){0.95}}
\put(2.75,.75){$F(\gamma)$}
\put(1.45,.75){$F(\pi)$}
\put(1.55,0.05){$\lambda$}
\put(2.75,0.05){$F(\alpha)$}

 \put(0.85,.6){\vector(0,-1){0.452}}
 \put(2.1,.6){\vector(0,-1){.452}}
\put(3.52,.6){\vector(0,-1){0.452}}

\put(.92,0.315){$\theta_{F,G}$}
 \put(2.17,0.315){$1_{F(A)}$}
\put(3.57,0.315){$F(\beta)$}
\end{picture}
\end{center}
If $F$ is in $\DD$, then $F(\pi)$ is monic; and hence $\theta_{F,G}$ is monic. 
 In order to show that $\theta$ is natural, assume that $\eta:G_1\To G_2$ is a morphism in $\modd(\C)$. Then for $i=1,2$, there is exact sequences $\Hom_{\C}(-,B_i)\To \Hom_{\C}(-,A_1)\To G_i\To 0$ such that $A_i$ and $B_i$ are objects in $\C$. Using Yoneda lemma there is the following commutative diagram of abelian groups with exact rows

\begin{center}\setlength{\unitlength}{3cm}
\begin{picture}(4.5,1.1)

\put(0.83,0.95){$\theta_{F,G_1}$}
\put(1.2,0.94){\vector(-3,-1){0.45}}
\put(1.78,0.94){\vector(3,-1){0.45}}
\put(1.2,0.95){$F(v(G_2))$}
\put(1.3,0.55){$F(v(\eta))$}

\put(1.78,0.26){\vector(3,-1){0.45}}
\put(1.2,0.26){\vector(-3,-1){0.45}}
\put(0.83,0.28){$\theta_{F,G_2}$}
\put(1.2,0.3){$F(v(G_1))$}
\put(1.3,.895){\vector(0,-1){0.452}}

 \put(0,0){\makebox(0,0){$0$}}
 \put(0.69,0){\makebox(0,0){$\Hom_{\modd(\C)}(G_1,F))$}}
\put(2.4,0){\makebox(0,0){$F(A_1)$}}
 \put(4,0){\makebox(0,0){$F(B_1)$}}

\put(0.0485,0){\vector(1,0){0.15}}
\put(1.2,0){\vector(1,0){1}}
\put(2.65,0){\vector(1,0){1.15}}

\put(0,0.7){\makebox(0,0){$0$}}
 \put(0.69,0.7){\makebox(0,0){$\Hom_{\modd(\C)}(G_2,F)$}}
\put(2.4,0.7){\makebox(0,0){$F(A_2)$}}
 \put(4,0.7){\makebox(0,0){$F(B_2)$}}

\put(0.05,0.7){\vector(1,0){0.15}}
\put(1.2,0.7){\vector(1,0){1}}
\put(2.65,0.7){\vector(1,0){1.15}}

 \put(0.5,.6){\vector(0,-1){0.452}}
 \put(2.45,.6){\vector(0,-1){.452}}
\put(4,.6){\vector(0,-1){0.452}}

\put(-0.12,.315){$\Hom_{\C}(\eta,F)$}
 \end{picture}
\end{center} 
The naturality $\theta_{F,G}$ in $F$ is similar.
\end{proof}

\medskip

\begin{Corollary}\label{c1}
Let $F$ and $G$ be in $\modd(\C)$. Then $F$ is in $\DD$ if and only if $\theta_{F,G}$ is monic for all functors $G$ in $\omodd(\C)$.
\end{Corollary}
\begin{proof}
Assume that $\theta_{F,G}$ is monic for all $G\in\omodd(\C)$ and that $B\To C\To 0$ is an exact sequence of objects of $\C$. Then there exists an exact sequence of functors $\Hom_{\C}(-,B)\To \Hom_{\C}(-,C)\To G\To 0$ such that $G\in\omodd(\C)$. Using Yoneda lemma, there is an exact sequence of abelian gropus $0\To \Hom_{\modd(\C)}(G,F)\To F(C)\To F(B)\To F(A)$. Since $\theta_{F,G}$ is isomorphism, $\Hom_{\C}(G,F)\cong F(v(G))=0$ and so $F$ is in $\DD$. The converse follows from the previous theorem.
\end{proof}

\medskip
\begin{Corollary}\label{c2}
Let $F\in \modd(\C)$. Then $F$ is half exact if and only if $\theta_{F,G}$ is epic for any $G\in\modd(\C)$. In particular, if $\omodd(\C)$ is Serre, then $F$ is left exact if and only if $\theta_{F,G}$ is isomorphism for all $G\in\modd(\C)$.
\end{Corollary}
\begin{proof}
Given an exact sequence $A\To B\To C\To 0$ such that $A\To B$ is a pseudo-kernel of $B\To C$, we have an exact sequence of functors $\Hom_{\C}(-,A)\To\Hom_{\C}(-,B)\To \Hom_{\C}(-,C)$. Putting $H=\Coker(\Hom_{\C}(-,A)\To\Hom_{\C}(-,B))$, since the functor $v$ preserves cokernels we have $v(H)\cong C$ and so $F(C)\To F(B)\To F(A)$ is an exact sequence of ablian groups if and only if $F(v(H))\To F(B)\To F(A)$ is exact. Now, we have the following commutative diagram with the bottom row exact  \\

\begin{center}\setlength{\unitlength}{3cm}
\begin{picture}(4.5,.7)
\put(0,0){\makebox(0,0){$0$}}
 \put(0.85,0){\makebox(0,0){$\Hom_{\C}(H,F))$}}
\put(2.1,0){\makebox(0,0){$\Hom_{\C}((-,B),F))$}}
 \put(3.5,0){\makebox(0,0){$\Hom_{\C}((-,A),F)).$}}

\put(0.1,0){\vector(1,0){0.297}}
\put(1.2473,0){\vector(1,0){0.33}}
 \put(2.58,0){\vector(1,0){0.45}}

 \put(0.9,0.7){\makebox(0,0){$F(v(H))$}}
\put(2.1,0.7){\makebox(0,0){$F(B)$}}
 \put(3.55,0.7){\makebox(0,0){$F(A)$}}

\put(1.2,.7){\vector(1,0){0.7}}
 \put(2.4,.7){\vector(1,0){0.95}}

 \put(0.85,.6){\vector(0,-1){0.452}}
 \put(2.1,.6){\vector(0,-1){.452}}
\put(3.52,.6){\vector(0,-1){0.452}}

\put(.92,0.315){$\theta_{F,H}$}
 \put(2.17,0.315){$\theta_{F,\Hom_{\C}(-,B)}$}
\put(1.99,0.315){$\cong $}
\put(3.57,0.315){$\theta_{F,\Hom_{\C}(-,A)}$}
\put(3.4,0.315){$\cong$}
\end{picture}
\end{center}
It is straightforward to prove that $\theta_{F,H}$ is epic if and only if the top row is exact. The second assertion follows by the first and the previous corollary, \cref{1} and \cref{leftex}.
\end{proof}

\medskip
In the rest of this section we assume that $\omodd(\C)$ is a Serre subcategry of $\modd(\C)$.
\begin{Corollary}
Let $F,G$ be in $\modd(\C)$. Then there is an isomorphism $\Ker\theta_{F,G}\cong \Ker \theta_{r(F),G)}$.
\end{Corollary}
\begin{proof}
There is an exact sequence of functors $0\To r(F)\To F\stackrel{\eta_F}\To h_{\C}(v(F))$. Then for any functor $G$ we have the following commutative diagram with exact rows 

\begin{center}\setlength{\unitlength}{3cm}
\begin{picture}(5.3,.5)

\put(.45,0.4){\makebox(0,0){$0$}}
\put(1.15,0.4){\makebox(0,0){$r(F)(v(G))$}}
 \put(2.6,0.4){\makebox(0,0){$F(v(G))$}}
\put(3.75,0.4){\makebox(0,0){$h_{\C}(v(F))(v(G)))$}}

\put(.45,0){\makebox(0,0){$0$}}
\put(1.15,0){\makebox(0,0){$\Hom_{\C}(G,r(F))$}}
 \put(2.6,0){\makebox(0,0){$\Hom_{\C}(G,F)$}}
\put(3.75,0){\makebox(0,0){$\Hom_{\C}(G,h_{\C}(v(F)))$}}
\put(3.815,0.2){\makebox(0,0){$\epsilon$}}

 \put(.556,.0){\vector(1,0){0.2}}
 \put(1.6,.0){\vector(1,0){0.6}}
\put(2.94,.0){\vector(1,0){.3}}

\put(.556,.40){\vector(1,0){0.25}}
 \put(1.5,.4){\vector(1,0){0.8}}
\put(2.86,.4){\vector(1,0){.382182}}

\put(1.1,.34){\vector(0,-1){.23}}
\put(1.35,0.2){\makebox(0,0){$\theta_{r(F),G}$}}
\put(2.52,.34){\vector(0,-1){0.23}}
\put(2.7,0.2){\makebox(0,0){$\theta_{F,G}$}}
\put(3.7,.3){\vector(0,-1){0.2}}
\end{picture}
\end{center}
The previous corollary implies that $\epsilon$ is isomorphism and so the result follows. 
\end{proof}
	
	\medskip
	
	\begin{Corollary}
	Let $F,G$ be in $\modd(\C)$. Then there is an isomorphism $\Ker\theta_{F,G}\cong \Ker \theta_{F,G/r(G)}$ and $\Coker\theta_{F,G/r(G)}$ is isomorphic to a subobject of $\Coker\theta_{F,G}$. In particular, if $r(F)=0$, then $\Coker\theta_{F,G/r(G)}\cong\Coker\theta_{F,G}$.
	\end{Corollary}
	\begin{proof}
	There are the following commutative diagram
	\begin{center}\setlength{\unitlength}{3cm}
\begin{picture}(5.3,.5)

\put(1.25,0.4){\makebox(0,0){$F(v(G/r(G)))$}}
 \put(2.7,0.4){\makebox(0,0){$F(v(G))$}}

\put(.45,0){\makebox(0,0){$0$}}
\put(1.25,0){\makebox(0,0){$\Hom_{\C}(G/r(G),F))$}}
 \put(2.7,0){\makebox(0,0){$\Hom_{\C}(G,F)$}}
\put(3.75,0){\makebox(0,0){$\Hom_{\C}(r(G),F)$}}

 \put(.556,.0){\vector(1,0){0.2}}
\put(2.05,0.085){\makebox(0,0){$(\pi,F)$}}
 \put(1.76,.0){\vector(1,0){0.6}}
\put(3.03,.0){\vector(1,0){.3}}

\put(1.85,0.46){\makebox(0,0){$\cong$}}
 \put(1.65,.4){\vector(1,0){0.8}}
\put(1.85,0.33){\makebox(0,0){$Fv(\pi)$}}

\put(1.1,.34){\vector(0,-1){.23}}
\put(1.35,0.2){\makebox(0,0){$\theta_{F,G/r(G)}$}}
\put(2.52,.34){\vector(0,-1){0.23}}
\put(2.7,0.2){\makebox(0,0){$\theta_{F,G}$}}

\end{picture}
\end{center}
It is clear that $\Ker\theta_{F,G}\cong\theta_{F,G/r(G)}$. This also implies that $\im\theta_{F,G/r(G)}\cong\im \theta_{F,G}$. Then we have the following commutative diagram with exact rows 
\begin{center}\setlength{\unitlength}{3cm}
\begin{picture}(5.3,.5)

\put(.45,0.4){\makebox(0,0){$0$}}
\put(1.15,0.4){\makebox(0,0){$\im\theta_{F,G/r(G)}$}}
 \put(2.56,0.4){\makebox(0,0){$\Hom_{\C}(G/r(G),F)$}}
\put(3.6,0.4){\makebox(0,0){$\Coker\theta_{F,G/r(G)}$}}
\put(4.3,0.4){\makebox(0,0){$0$}}

\put(.45,0){\makebox(0,0){$0$}}
\put(1.15,0){\makebox(0,0){$\im\theta_{F,G}$}}
 \put(2.6,0){\makebox(0,0){$\Hom_{\C}(G,F)$}}
\put(3.65,0){\makebox(0,0){$\Coker\theta_{F,G}$}}
\put(4.3,0){\makebox(0,0){$0$}}
\put(3.9,0.2){\makebox(0,0){$\overline{(\pi,F)}$}}

 \put(.556,.0){\vector(1,0){0.35}}
 \put(1.38,.0){\vector(1,0){0.86}}
\put(2.97,.0){\vector(1,0){.38}}
\put(3.97,.0){\vector(1,0){.273}}

\put(.556,.40){\vector(1,0){0.25}}
 \put(1.5,.4){\vector(1,0){0.56}}
\put(3.03,.4){\vector(1,0){.17}}
\put(4,.4){\vector(1,0){.23}}

\put(1.1,.34){\vector(0,-1){.23}}

\put(2.52,.34){\vector(0,-1){0.23}}
\put(2.7,0.2){\makebox(0,0){$(\pi,F)$}}
\put(3.7,.3){\vector(0,-1){0.2}}
\end{picture}
\end{center}
which implies that $\Ker\overline{(\pi,F)}=0$. If $r(F)=0$, then $\Hom_{\modd(\C)}(r(G),F)=0$ and then $(\pi,F)$ is isomorphism so that $\Coker\overline{(\pi,F)}=0$. 
	\end{proof}

In the rest of this paper we study the injective objects in $\modd(\C)$. 
\medskip
\begin{Lemma}\label{inje}
If $G$ is a functor in $\modd(\C)$ such that $\Ext_{\modd(\C)}^i(F,G)=0$ for $i=1,2$ and all $F\in\omodd(\C)$, then $G$ is injective in $\modd(\C)$.
\end{Lemma}
\begin{proof}
For any functor $F$ in $\modd(C)$ consider the exact sequence of functors $$0\to r(F)\To F\stackrel{\eta_F}\To h_{\C}v(F)\To \Coker(\eta_F)\To 0$$
in which $\Ker(\eta_F)$ and $\Coker(\eta_F)$ are in $\omodd(\C)$. We observe that $h_{\C}v(F)=\Hom_{\C}(-,v(F))$ is projective in $\modd(\C)$ and hence $\Ext_{\modd(\C)}^1(h_{\C}v(F),G)=0$. Then applying the functor $\Hom_{\modd(\C)}(-,G)$ to the above exact sequence and using the assumption, we deduce that $\Ext_{\modd(\C)}^1(F,G)=0$. 
\end{proof}
\medskip
\begin{Proposition}\label{rinj}
If $G\in\modd(\C)$ is injective, then so is $r(G)$.
\end{Proposition}
\begin{proof}
By \cref{adje}, since $r:\modd(\C)\To\omodd(\C)$ is a right adjoint functor of the inclusion functor $i:\omodd(\C)\To\modd(\C)$, $r(G)$ is injective in $\omodd(\C)$; and hence $\Ext_{\modd(\C)}^i(F,r(G))=0$ for all $i>0$ and all $F\in\omodd(\C))$. Now, \cref{inje} implies that $r(G)$ is injective. 
\end{proof}

\medskip

\begin{Proposition}\label{inje}
Let $\C$ be abelian. If $A$ is an object of $\C$ such that ${\rm id}A=n$, then $\Ext^n_{\C}(-,A)$ is injective in $\modd(\C)$. Moreover, if $n\geq 1$, then $\Ext^n_{\C}(-,A)$ is in $\omodd(\C)$. 
\end{Proposition}
\begin{proof}
The case $n=0$, since by \cref{ab}, $v$ is left exact, it is straightforward to show that $h_{\C}(A)=\Hom_{\C}(-,A)$ is injective in $\modd(\C)$. Assume that $n=1$ and consider an injective resolution $0\To A\To E_0\To E_1\To 0$ for $A$. By virtue of the case $n=0$, there exists an exact sequence of functors $$0\To\Hom_{\C}(-,A)\To \Hom_{\C}(-,E_0)\To \Hom_{\C}(-,E_1)\To \Ext_{\C}^1(-,A)\To 0$$ which implies that $\Ext_{\C}^1(-,A)\in\omodd(\C)$. Assume that $G=\Coker(\Hom_{\C}(-,A)\To \Hom_{\C}(-,E_0))$. For any functor $F\in\C$, using the adjointness $\Hom_{\C}(-,E_i)$ are injective for $i=0,1$, we conclude that $\Ext_{\modd(\C)}^1(F,\Ext_{\C}^1(-,A))=\Ext_{\modd(\C)}^3(F,(-,A))=0$ because pd$F\leq 2$. Thus
$\Ext_{\C}^1(-,A)$ is injective. Assume that $n>1$ and hence there exists an exact sequence of objects $0\To A\To E\To B\To 0$ in $\C$ such that $E$ is injective and id$B=n-1$. An easy induction and the isomorphism $\Ext_{\C}^n(-,A)\cong \Ext_{\C}^{n-1}(-,B)$ imply that 
 $\Ext_{\C}^n(-,A)$ is in $\omodd(\C)$ and it is injective in $\C$.
\end{proof}
\medskip

		\begin{Remark}
			We notice that if $\CA$ is an abelian category with enough projective, then it follows from [AB, Proposition 1.8] that $h_{\C}v(F)\cong R^0F$ and so $\Ker(F\stackrel{\eta_F}\To h_{\C}v(F))=r(F)$ is projective stable. We also observe that by using [CE, Chap V. Theorem 5.3] that for any $n>0$, we have an equivalence of functors $R^nF(-)\cong {\rm Ext}^n(-,v(F))$ where $R^n$ is the $n$-th right derived functor of $F$. We observe that $\underline{\modd}(\C)$ consists of all projective stable functors. It is clear that $F\in\omodd(\C)$ if and only if $R^0F=0$ and hence in this case $R^nF=0$ for all $n\geq 0$.
				
	If $0\To F_1\To F_2\To F_3\To 0$ is an exact sequence in $\modd(\C)$, then there is a long exact sequence of right derived 
	functors in $\Modd(\C)$. To be more precise, we have an exact sequence $0\To v(F_1)\To v(F_2)\To v(F_3)\To 0$ in $\C$ and so there is a long exact sequence in $\Modd(\C)$ $$0\To R^0(F_1)\To R^0(F_2)\To R^0(F_3)\To R^1(F_1)\To\dots. $$
			We observe that if $\C$ has enough injective objects, then $R^nF\in\modd(\C)$ for any $F\in\modd(\C)$. Because, there is an exact sequence $0\To v(F)\To E\To C\To 0$ such that $E$ is injective. Then there is an exact sequence $\Hom_{\C}(-,E)\To \Hom_{\C}(-,C)\To R^1F\To 0$ so that $R^1F\in\modd(\C)$. An induction on $n$ implies that $R^nF\in\modd(\C)$ for all $n\geq 1$. 			
	\end{Remark}		

\section{limits in coherent functors}
Throughout this section $\omodd(\C)$ is a Serre subcategory of $\modd(\C)$.

A non-empty category $\CI$ is said to be filtered provided that for each pair of objects $\lambda_1,\lambda_2$ of $\CI$, there are morphisms $\phi_i: \lambda_i\to \mu$ for some $\mu\in \CI$, and for each pair of morphisms $\phi_1,\phi_2: \lambda\to \mu$, there is a morphism $\psi:\mu\To v$ with $\psi\phi_1 = \psi\phi_2$. We call the colimit $\underrightarrow{\rm lim}X_{\lambda}$ of a functor $X: \CI \To \C, \lambda\mapsto X_{\lambda}$ a {\it direct limit} if $\CI$ is a skeletally small filtered category. We denote by $\overrightarrow{\omodd(\C)}$ the full subcategory of $\Modd(\C)$ which consists of direct limit $\underrightarrow{\rm lim} F_i$ with each $F_i\in\omodd(\C)$.

\medskip
\begin{Proposition}\label{td}
 There are the following equalities:
$$\{X\in\Modd(\C)|\Hom_{\Modd(\C)}(\omodd(\C),X)=0\}$$$$=\{X\in\Modd(\C)|\Hom_{\Modd(\C)}(\overrightarrow{\omodd(\C)},X)=0\}=\overrightarrow{\DD}.$$
\end{Proposition}
\begin{proof}
The first equality is clear. Clearly $\overrightarrow{\DD}\subseteq\{X\in\Modd(\C)|\Hom_{\Modd(\C)}(\omodd(\C),X)=0\}$. Conversely, assume that $X\in\{X\in\Modd(\C)|\Hom_{\Modd(\C)}(\omodd(\C),X)=0\}$. By [Cr, Lemma 4.1], it suffices to show that for any $Y\in\modd(\C)$, any morphism $Y\To X$ factors through some $D\in\DD$. Consider the exact sequence $0\To r(Y)\To Y\To Y/r(Y)\To 0$ where $r(Y)\in\omodd(\C)$ and $Y/r(Y)\in\DD$. Since $\Hom_{\Modd(\C)}(r(Y),X)=0$, the morphism $Y\To X$ factors through $Y\To Y/r(Y)$. 
\end{proof}

\medskip
Let $\CA$ be an abelian category. A torsion theory $(\cT,\cF)$ of $\CA$ is said to be of {\it finite type} if the corresponding right adjoint $r:\CA\To\cT$ of the inclusion functor $\iota:\cT\To\CA$ commutes with direct limits. Let $\cS$ be a Serre subcategory of $\CA$ and let $\CA/\cS$ be the corresponding quotient category with the canonical functor $q:\CA\to \CA/\cS$ (for details we refer the readers to [G, P]). The subcategory $\cS$ of $\CA$ is called {\it localizing} provided that $q$ admits a right adjoint functor $s:\CA/\cS\To \CA$ which is called {\it section functor}. For each $X\in\CA$, the adjointness induces a natural morphism $\psi_X:X\To s\circ q(X)$ with $\Ker\psi_X,\Coker\psi_X\in\cS$. An object $X\in\CA$ is said to be $\cS$-{\it closed} provided $\psi_X$ is an isomorphism.

\medskip
\begin{Proposition}\label{loca}
 $\overrightarrow{\omodd(\C)}$ is a localizing subcategory of finite type of $\Modd(\C)$. 
\end{Proposition}
\begin{proof}
As $\Modd(\C)$ is locally coherent, the result follows by [K3, Theorem 2.8].
\end{proof}

	\medskip
	 
	Given a subclass $\cS$ of an abelian category $\CA$, we recall that the {\it right prependicular category} $\cS^\bot$ of $\cS$ is the full subcategory of $\CA$ consisting of all objects $M$ satisfying $\Hom_{\CA}(X,M)=\Ext^1_{\CA}(X,M)=0$ for all $X\in\cS$. If $\cS$ is a localizing subcategory of $\CA$, by virtue of [G, III.2, Lemme 1], the right prependicular category $\cS^\bot$ coincides with the full subcategory of $\cS$-closed objects.   
	
		\begin{Proposition}\label{st}
		There are the equalities $\omodd(\C)^\bot=\overrightarrow{\omodd(\C)}^\bot=\overrightarrow{\cL}$ where $\cL$ is a subclass of $\modd(\C)$ consisting of left exact functors. 
		\end{Proposition}
		\begin{proof}
 The first equality follows from [K3, Corollary 2.11]. According to \cref{loca}, $\overrightarrow{\omodd(\C)}$ is a localizing subcategory of $\Modd(\C)$ and so by [G, III.2, Lemme 1], it suffices to show that $X$ is $\overrightarrow{\omodd(\C)}$-closed. Assume that $q:\Modd(\C)\To \Modd(\C)/\overrightarrow{\omodd(\C)}$ is the quotient functor with the section functor $s:\Modd(\C)/\overrightarrow{\omodd(\C)}\To \Modd(\C)$. Since, by \cref{loca}, $\overrightarrow{\omodd(\C)}$ is of finite type, it follows from [K3, Lemma 2.4] that $s$ preserves direct limits. If $X\in\overrightarrow{\cL}$, then $X=\underrightarrow{\rm lim} X_i$ where $X_i\in\cL$ for each $i$. Using \cref{leftex}, we have $\Ext_{\Modd(\C)}^i(\omodd(\C),X_i)=0$ for $i=0,1$. For each $i$, we have an exact sequence $0\To \Ker\psi_{X_i}\To X_i\stackrel{\psi_{X_i}}\To s\circ q(X_i)\To \Coker\psi_{X_i}\To 0$ where $\Ker\psi_{X_i}$ and $\Coker\psi_{X_i}$ are in $\overrightarrow{\omodd(\C)}$. Now the previous argument  implies that $\Ker\psi_{X_i}=0$. If $\Coker\psi_{X_i}=\underrightarrow{\rm lim} C_j$ with $\C_j\in\omodd(\C)$, for each $j$, applying the functor $\Hom_{\Modd(\C)}(C_j,-)$ to the above exact sequence and the fact that $\Ext^1(C_j,X_i)=0$ we conclude that each the canonical morphism $\phi_j:C_j\To\Coker\psi_{X_i}$ factors through $s\circ q(X_i)$ so that $\phi_j=0$. Thus $\Coker\psi_{X_i}=0$ and so $\psi_{X_i}$ is isomorphism for each $i$. Now since $s$ and $q$ preserves direct limits, we have $X\stackrel{\psi_X}\cong s\circ q(X)$. Conversely assume that $X\in\omodd(\C)^\bot$ and $Y\in\modd(\C)$. Then applying the functor $\Hom_{\Modd(\C)}(-,X)$ to the exact sequence $0\To r(Y)\To Y\stackrel{\eta_Y}\To h_{\C}v(Y)\To \Coker\eta_Y\To 0$, we deduce that any morphism $Y\To X$ factors through $\eta_Y:Y\To h_{\C}v(Y)$. Since $h_{\C}v(Y)\in\cL$, the result follows from [Cr, Lemma 4.1].
		\end{proof}
		\medskip
		Let $\cM$ be a subcategory of $\Modd(\C)$ and let $F$ be a functor in $\Modd(\C)$. A morphism $f:F\To M$ with $M\in\cM$ is an $\cM$-{\it preenvelope} (or a left $\cM$-approximation) of $F$ provided that the induced abelian group homomorphism $\Hom_{\Modd(\C)}(f,X):\Hom_{\Modd(\C)}(M,X)\to\Hom_{\Modd(\C)}(F,X)$ is surjective for each $X\in\cM$.

		\begin{Corollary}\label{t}
	If $\C$ has direct limits, then any functor in $\overrightarrow{\cL}$ has an $\cL$-preenvelope.
	\end{Corollary}
\begin{proof}
Assume that $F\in\overrightarrow{\cL}$. Then $F=\underrightarrow{\rm lim} F_i$, where $F_i\in{\cL}$. Using \cref{2}, there exists some object $A_i$ such that $F_i=\Hom_{\C}(-,A_i)$ for each $i$. Assume that $\sigma_i:F_i\To F$ is the canonical morphism for each $i$. If $\bF:\CI\to \Modd(\C)$ defined as $\bF(i)=F_i$, is the corresponding functor where $\CI$ is an skeletally small filtered category, then $v(\bF):\CI\To \C$ defined as $v(\bF)(i)=A_i$ is a functor and since $\C$ has direct limits, $\underrightarrow{\rm lim}v(\bF)=\underrightarrow{\rm lim} A_i$ exists with the canonical morphism $\phi_i:A_i\to\underrightarrow{\rm lim} A_i$ for each $i$. We observe that there exists a unique morphism $\delta:F\To \Hom_{\C}(-,\underrightarrow{\rm lim} A_i)$ such that $\delta\circ\sigma_i= h_{\C}(\phi_i)$. We assert that $\delta$ is an $\cL$-preenvelope. Suppose that $G=\Hom_{\C}(-,C)$ is in $\cL$ and $f:F\To G$ is a morphism. Then there exists a unique morphsim $\epsilon:\underrightarrow{\rm lim} A_i\To C$ such that $\epsilon\circ\phi_i=v(f\circ\sigma_i)$. Thus for each $i$, we have $h_{\C}(\epsilon)\circ\delta\circ\sigma_i=h_{\C}(\epsilon)\circ h_{\C}(\phi_i)=h_{\C}(v(f\circ\sigma_i))=f\circ\sigma_i$. The unicity implies that $h_{\C}(\epsilon)\circ\delta=f$. 
\end{proof}
\medskip
The following results shows that $\overrightarrow{\cL}$ is the class of right exact functors in $\Modd(\C)$. 
  
\begin{Proposition}\label{lefex}
 Let $\C$ has pseudo-kernels and $F\in\Modd(\C)$. Then $F\in\overrightarrow{\cL}$ if and only if $F$ is left exact. 
\end{Proposition}
\begin{proof}
If $F\in\overrightarrow{\cL}$, then $F=\underrightarrow{\rm lim} F_i$, where $F_i\in{\cL}$. Then $F$ is left exact as the direct limits are exact in the categorey of abelian groups. Conversely assume that $F$ is left exact. For any $G\in\omodd(\C)$, using \cref{epimono} there exists an epimorphism $B\To A$ in $\C$ such that $\Hom_{\C}(-,B)\To\Hom_{\C}(-,A)\To G\To 0$ is exact. Assume that $C\To B$ is a pseudo-kernel of $B\To A$. It follows from \cref{copse} that $C\To B\To A\To 0$ is exact and since $F$ is left exact, $0\To F(B)\To F(A)\To F(C)$ is an exact sequence of abelian groups. On the other hand, there is an exact sequence $\Hom_{\C}(-,C)\To \Hom_{\C}(-,B)\To\Hom_{\C}(-,A)\To G\To 0$. Applying $\Hom_{\Modd(\C)}(-,F)$ to this exact sequence and using Yoneda lemma, we deduce thtat $\Hom_{\Modd(\C)}(G,F)=\Ext^1_{\Modd(\C)}(G,F)=0$ so that $F\in\omodd(\C)^\bot$. Now the result follows by using \cref{st}.  
\end{proof}

\medskip
\begin{Proposition}\label{ex}
For any functor $F\in\Modd(\C)$, there exists a morphism $\eta_F:F\To L(F)$ such that $L(F)\in\overrightarrow{\cL}$ and $\Ker\eta_F$ and $\Coker\eta_F$ are in $\overrightarrow{\omodd(\C)}$.   
\end{Proposition}
\begin{proof}
Since $\Modd(\C)$ is locally coherent, we have $F=\underrightarrow{\rm lim} F_i$ where $F_i\in\modd(\C)$ for each $i$. Then for each $i$, there is an exact sequence $$0\To r(F_i)\To F_i\stackrel{\eta_{F_i}}\To \Hom_{\C}(-,v(F_i))\To c(F_i)\To 0$$
where $r(F_i), c(F_i)\in\omodd(\C)$ for all $i$. If $\CI$ is an skeletally small filtered category such that $\bF:\CI\To\Modd(\C)$ is a functor with $\bF(i)=F_i$, then $r(\bF):\CI\to\Modd(\C)$ defined as $r(\bF)(i)=r(F_i)$ is a functor which has direct limit $\underrightarrow {\rm lim} r(F_i)$ as $\Modd(\C)$ is a Grothendeick categorey by [B, Satz 1.5]. The same holds for $\underrightarrow{\rm lim} c(F_i)$. Then we have an exact sequence of functors 
$$0\To \underrightarrow{\rm lim} r(F_i)\To F\stackrel{\eta_{F}}\To\underrightarrow{\rm lim} \Hom_{\C}(-,v(F_i))\To\underrightarrow{\rm lim} c(F_i)\To 0$$
in which $L(F)=\underrightarrow{\rm lim} \Hom_{\C}(-,v(F_i))\in\overrightarrow{\cL}$.
\end{proof}

We consider a functor $r:\Modd(\C)\To\overrightarrow{\omodd(\C)}$ defined as $r(F)=\Ker\eta_F$ for any functor $F\in\Modd(\C)$. It follows from \cref{st} that $r(F)$ is the largest subfuncor of $F$ contained in $\overrightarrow{\omodd(\C)}^\bot$. Since $\overrightarrow{\omodd(\C)}$ is a localizing subcategory of $\Modd(\C)$, the functor $r$ is idempotent and radical, and hence $r:\Modd(\C)\To \Modd(\C)$ is left exact. Using again \cref{st}, the functor $r:\Modd(\C)\To \overrightarrow{\omodd(\C)}$ is a right adjoint functor of the inclusion functor. 
\medskip

\begin{Corollary}\label{torm}
$(\overrightarrow{\omodd(\C)},\overrightarrow{\DD})$ is a hereditary torsion theory of $\Modd(\C)$.
\end{Corollary}
\begin{proof}
Assume that $\Hom_{\Modd(\C)}(X,Y)=0$ for all $Y\in\overrightarrow{\DD}$. It is clear that $X/r(X)\in\overrightarrow{\DD}$ and so $X\in\overrightarrow{\omodd(\C)}$. The other part follows by \cref{td}.
\end{proof}

\medskip

\begin{Proposition}\label{adj}
There exists a right exact functor $v^{*}:\Modd(\C)\To \C$ such that $v^*|_{\modd(\C)}=v$ and it is a left adjoint functor of $h_{\C}:\C\to\Modd(\C)$.
\end{Proposition}
\begin{proof}
According to [K2, Universal Property 5.6] the functor $v^*$ exists such that commutes with direct limits so that it is right exact. For the second assertion, for any $F\in\Modd(\C)$ and $A\in\C$, since $\Modd(\C)$ is locally coherent, we have $F=\underrightarrow{\rm lim} F_i$, where $F_i\in{\modd(\C)}$; and hence there are the following isomorphisms 
$$\Hom_{\C}(v^*(F),A)=\Hom_{\C}(\underrightarrow{\rm lim}v(F_i),A)\cong \underleftarrow{\rm lim}\Hom_{\C}(v(F_i),A)$$$$\cong \underleftarrow{\rm lim}\Hom_{\modd(\C)}(F_i,h_{\C}(A))\cong \Hom_{\modd(\C)}(F,h_{\C}(A)).$$
\end{proof}

\medskip

\begin{Proposition}
If $\C$ has direct limits, then every $F\in\Modd(\C)$ has $\cL$-preenvelopes.
\end{Proposition}
\begin{proof}
Since $\Modd(\C)$ is locally coherent, we have $F=\underrightarrow{\rm lim} F_i$ where $F_i\in\modd(\C)$ for each $i$. By the proof of \cref{ex}, we have $L(F)=\underrightarrow{\rm lim} \Hom_{\C}(-,v(F_i))$ and using \cref{t}, the morphism $\delta:\underrightarrow{\rm lim} \Hom_{\C}(-,v(F_i))\To h_{\C}v^*(F)$ is an $\cL$-preenvelope. We claim that $\delta\circ\eta_F:F\To h_{\C}v^*(F)$ is an $\cL$-preenvelope. For any $G\in\cL$, applying the functor $\Hom_{\Modd(\C)}(-,G)$ to the exact sequence $$0\To \underrightarrow{\rm lim} r(F_i)\To F\stackrel{\eta_{F}}\To L(F)\To\underrightarrow{\rm lim} c(F_i)\To 0$$ and using \cref{st}, we deduce that for any morphism $\theta:F\To G$, there exists a morphism $\epsilon:L(F)\To G$ such that $\epsilon\circ\eta_F=\theta$ and using \cref{t}, there exists a morphism $\alpha:h_{\C}v^*(F)\To G$ such that $\alpha\circ\delta=\epsilon$. Thus $\alpha\circ\delta\circ\eta_F=\theta$.   
\end{proof}



\end{document}